\magnification \magstep1 \openup 2\jot \def \qed {\vrule height6pt 
width6pt depth0pt}
\def \sgn {{\rm sgn\,}}
\def \rank {{\rm rank\,}}
\def \dim {{\rm dim\,}}
\centerline{\bf{Factorizations of natural embeddings of $l_p^n$ into $L_r$, 
II}}

\centerline { by}

\centerline { T.~Figiel, W.~B.~Johnson\footnote {*}{ Supported in part by NSF 
DMS-87-03815} and G.~Schechtman\footnote { **}{ Supported in part by the 
U.S.-Israel BSF} }

{\ } 

\noindent{\bf Introduction}

In this continuation of [FJS], we show that in some situations considered 
in [FJS], conclusions of certain theorems can be strengthened.
More explicitly, suppose that $T$ is an operator from some Banach space 
 into $L_1$ which factors through some $L_1$-space $Z$ as $uw$ and 
normalized so that  $\Vert w \Vert = 1$.  In Corollary 12.A we show that 
if $T$ is the inclusion mapping from a ``natural" $n$-dimensional Hilbertian 
subspace  of $L_1$ into  $L_1$, then $u$ well-preserves a copy of 
$l_1^k$ with $k$ exponential in $n$ (where ``well" and the base of the 
exponent depend on $\Vert u \Vert$  and on a quantitative measure of 
``naturalness").  This improves the result of [FJS] that the same hypotheses 
yield that $l_1^k$ well-embeds into $uZ$.  Corollary 12.B gives a similar 
improvement of Corollary 1.5 in [FJS]; that is, Corollary 12.B is the same 
as Corollary 12.A except that the operator $T$ is assumed to be a  
mapping from a space  whose dual has controlled cotype into $L_1$ 
which acts like a quotient mapping relative to a ``natural" Hilbertian 
subspace of $L_1$.  

Corollary 20 strengthens the conclusion of Proposition 4.3 in [FJS] in a 
similar manner; it states that an operator from a $C(K)$ space which 
well-preserves a copy of $l_2^n$ also well-preserves a copy of $l_{\infty}^k$ 
with $k$ exponential in $n$ (rather than just have rank which is exponential 
in $n$). This can be 
viewed as a finite dimensional analogue of a particular case of a result of 
Pe\l czy\'nski [Pe1] stating that every non weakly compact operator from 
a $C(K)$ space preserves a copy of $c_0$.

In Theorem 21 we apply the earlier results in order to prove that for 
each $m$ there is an $m$-dimensional normed space $G$ such that any 
superspace of $G$ with a good unconditional basis must contain a copy of 
$l_\infty ^k$ with $k$ exponential in $m$.

We thank J. Bourgain for pointing us in the right direction on the material 
presented here.  After proving the results in [FJS], we suggested to him  
that there might be a translation invariant operator $T$ of bounded norm 
on $L_1$ of the circle which is the identity on the span of the first $n$ 
Rademacher functions and which does not preserve $l_1^k$ with $k$ 
exponential in $n$.  By disproving this conjecture, Bourgain started us 
thinking that Corollary 12.A was true.

{\ } 

\noindent{\bf Preliminaries, a quantitative version of Rosenthal's lemma}

In this section we prove, in Proposition 1 below, a quantitative version of 
Maurey's formulation [M] of Rosenthal's lemma [R] stating that an operator 
into  $L_1(\mu)$  either factors through an  $L_p(\nu)$   space for some   
$p>1$  via a change of density or is of type no better than 1.  Our approach
is in fact close in spirit to Rosenthal's original argument which (unlike some
later arguments) was basically quantitative in nature. Recall that for 
$1\le p<q\le\infty$  and for an operator    $u\colon Z\rightarrow L_p(\mu)$  
$$C_{p,q}(u)=\inf\{ \Vert h\Vert _s\Vert h^{-1}u\colon Z\rightarrow 
L_q\Vert \}$$ where the {\it inf}\/ is over all changes of measure  $h$; i.e., 
over all   $0<h\in L_s(\mu)$   where  ${1\over q}+{1\over s}={1\over p}$. 

In the statement of Proposition 1 and elsewhere in this paper $t^*$ 
denotes $t\over{t-1}$.

\proclaim Proposition 1. {\sl Let  $1<p<q\le\infty$ and $\sigma=1-{q^*\over 
p^*}$. Let 
$T\colon Z\rightarrow L_1(\mu)$, where 
$\mu$  is  a probability measure. If  $ C_{1,p}(T)=K$ 
and   $\Vert T \colon Z\rightarrow L_p(\mu) \Vert = {\cal C}K$, then for 
some  
$m>\Bigl({K\over {2^{2+{1\over p}}
{\cal C}\Vert T\Vert }}\Bigr)^{p^*}$ there exist $z_1,\dots,z_m$ in the unit 
ball of $Z$ and mutually disjoint measurable sets $F_1,\dots,F_m$ such that 
for $i=1,\dots,m$ one has}
$$\Vert 1_{F_i}\,Tz_i\Vert_1^{\sigma} \Vert 1_{F_i}\,Tz_i\Vert_q^{1-\sigma} 
\ge 2^{-
({1\over p} + 1)} K.\eqno(1) $$

For the proof we need two lemmas.

\proclaim Lemma 2. {\sl Let  $g \in L_1$ with $\Vert g \Vert \le 1$. 
Suppose $E$ is a 
$\mu$-measurable set, $1<p<q\le\infty$ and} 
$$\Vert 
1_{\lower3pt\hbox{$\scriptstyle \sim E$}}\,g   \Vert_p > \kappa >0. $$ 
{\sl Then there exists a measurable set $F$, 
$F\cap E =\emptyset$, such that} $$\mu(F)<\Bigl({2^{1/p}\over 
\kappa}\Bigr)^{p^*}, $$ 
$$\Vert 1_{F}\,g \Vert_1^{\sigma} \Vert 1_{F}\,g \Vert_q^{1-\sigma} > 
2^{-{1\over p}} \kappa. $$

   {\noindent \bf Proof.} Without loss of generality we can assume that $g 
\ge 0$. Set 
$F=[g>\gamma] \sim E$, where $\gamma >0$ is 
defined below.  Observe that 
$$\int_{\sim F}g^p\le \gamma^{p-1}\int_{\sim F}g\le \gamma^{p-1} $$ 
and hence, by H\"older's inequality, 
$$\kappa^p -\gamma^{p-1}<\int_Fg^p\le \Bigl(\int_F 
g\Bigr)^{1-t} \Bigl(\int_F g^q\Bigr)^t, $$ where $t=(p-1)/(q-1)$. Since 
$\mu(F)< 1/\gamma$, we can fulfill both conditions of the lemma by 
choosing $\gamma$ so that $\gamma^{p-1}={1\over 2}\kappa^p$.\hfill\qed

  \proclaim Lemma 3. {\sl Let $T\colon Z\rightarrow L_1(\mu)$, $\mu$ 
being a 
probability measure. Suppose $1<p<\infty$ and} $$ C_{1,p}(T)=K,\quad\quad 
\Vert T \colon Z\rightarrow L_p(\mu) \Vert = {\cal C}K. $$ {\sl If 
$0<\kappa<K$, 
$\eta\ge 0$ and ${\cal C}\eta^{1/{p^*}}\le 1-\kappa/K$, then 
$\mu(E)\le\eta$ implies} that there exists $z\in Ball(Z)$ such that 
$\Vert 
1_{\lower3pt\hbox{$\scriptstyle \sim E$}}\,Tz   \Vert_p > \kappa$.

   {\noindent \bf Proof.} Without loss of generality we may assume that 
$\mu(E)>0$. Observe that for any measurable set $A$ one has $$ 
C_{1,p}(1_{\lower3pt\hbox{$\scriptstyle A$}}\,T\colon Z\rightarrow 
L_1(\mu))\le \mu(A)^{1/p^*}\Vert 
1_{\lower3pt\hbox{$\scriptstyle A$}}\,T 
\colon Z\rightarrow L_p(\mu) \Vert. $$ The lemma follows by using this 
observation for $A=\sim E$ and $A=E$, because $$\Vert 
1_{{\lower3pt\hbox{$\scriptstyle \sim E$}}}\,T 
\colon Z\rightarrow L_p(\mu) \Vert > 
C_{1,p}(1_{{\lower3pt\hbox{$\scriptstyle \sim 
E$}}}\,T) \ge C_{1,p}(T) -
C_{1,p}(1_{\lower3pt\hbox{$\scriptstyle E$}}\,T) $$ 
$$\ge K-\mu(E)^{1/p^*}\Vert 1_{\lower3pt\hbox{$\scriptstyle E$}}\,T \colon 
Z\rightarrow 
L_p(\mu) \Vert \ge K(1-{\cal C}\eta^{1/p^*})\ge\kappa. \eqno{\qed}$$ 

{\noindent \bf Proof of Proposition 1.} Clearly, we may assume that \ $\Vert 
T\Vert =1$. Put \ $\kappa = {1\over 2}K$\ ,\break $\eta = (2{\cal C})^{-p^*}$, 
$\delta = \bigl( {{2^{1/p}} \over \kappa} \bigr)^{p^*}$  and let us start with 
$E=\emptyset$. Since 
$$\Vert 1_{\lower3pt\hbox{$\scriptstyle {\sim E}$}} \,T \colon Z\rightarrow 
L_p(\mu) \Vert > \kappa , $$ 
using Lemma 2 we can define $z_1$ and $F_1$ so that 
$\mu(F_1)<\delta$ and $\Vert 1_{F_1}\,Tz_1\Vert_1$ satisfies (1).   Suppose 
now that, for some $i\ge 1$, we have already defined $z_1,\dots,z_i$ and  
$F_1,\dots,F_i$. Let $E=\bigcup_{j\le i}F_j$. As long as $\mu(E)\le\eta$, 
Lemma 3 guarantees that we can use Lemma 2 again in order to choose 
$z_{i+1}$ and $F_{i+1}$ so that $F_{i+1}\cap E =\emptyset$, 
$\mu(F_{i+1})<\delta$ and $\Vert 1_{F_{i+1}}\,T z_{i+1} \Vert_1$ satisfies 
the estimate (1).   Therefore this procedure can be applied more than 
$\eta/\delta$ times. Since we have been assuming $\Vert T\Vert =1$, this 
yields the promised lower estimate for $m$ and completes the proof of the 
proposition.\hfill\qed

The next proposition shows that we can actually get a somewhat stronger 
conclusion to Proposition 1; namely, for some $k$ proportional to $m$, the 
identity on $l_1^k$ can be factored through T. Recall that, for a pair of linear 
operators  $T\colon X\rightarrow Y$ and  $U\colon X_1\rightarrow Y_1$  , 
the 
factorization constant of $U$ through $T$ is defined to be 
$$\gamma_{\lower3pt\hbox{$\scriptstyle T$}}(U) = \inf \{\Vert A\Vert 
\Vert B\Vert :  A\colon X_1\rightarrow X, \ B\colon Y\rightarrow Y_1, \ 
U=BTA  \}. $$ 
We let 
$\gamma_{\lower3pt\hbox{$\scriptstyle T$}}(U) = \infty$ if no such 
factorization exists. We also put $$\gamma_{\lower3pt\hbox{$\scriptstyle 
T$}}(Z) = \gamma_{\lower3pt\hbox{$\scriptstyle T$}}(id_Z\colon 
Z\rightarrow Z). 
$$

\proclaim Proposition 4. {\sl Let $T\colon Z\rightarrow L_1(\mu)$ be a 
bounded 
linear operator. Suppose that $z_1,\dots,z_m$ are in the unit ball of $Z$ and 
$F_1,\dots,F_m$ are mutually disjoint $\mu$-measurable sets such that, for 
$i=1,\dots,m$,} 
$$ \Vert 1_{\lower3pt\hbox{$\scriptstyle {F_i}$}}Tz_i\Vert\ge\delta\Vert 
T\Vert>0. $$ 
{\sl 
Then for some $k\ge {1\over 8}\delta m $  there exist linear operators 
$A\colon l_1^k\rightarrow Z$ and $B\colon L_1(\mu)\rightarrow l_1^k$ such 
that $BTA= 
id_{l_1^k}$ and $\Vert A\Vert\Vert B\Vert \Vert T\Vert\le 2\delta^{-1}$; 
i.e., $\gamma_{\lower3pt\hbox{$\scriptstyle T$}}(l_1^k)\le 2\delta ^{-
1}\Vert T\Vert ^{-1}$ for some $k\ge {1\over 8}\delta m$.}

 For the proof we need two basically known lemmas (see  [JS]). 

 \proclaim Lemma 5.  {\sl  Let $x_1,\dots,x_m$ be elements of $L_1(\mu)$ 
and let $A_1,\dots,A_m$ be mutually disjoint  $\mu$-measurable sets. If 
$1<k\le {1\over 2}m$, then there exists a subset $D\subset \{1,\dots,m \}$ 
with $\vert D\vert =k$ such that for each} $i\in D$ $$\sum_{j\in D\sim  
\{i\}}\int_{A_j}\vert x_i \vert\,d\mu\le (2k-1){m\choose 2}^{-
1}\sum_{i=1}^m\Vert x_i\Vert . $$ 

   {\noindent \bf Proof.} Setting $a_{ij}=\int_{A_j}\vert x_i\vert\,d\mu$ , for 
$i,j=1,\dots,m$, and $\alpha = \sum_{1\le i\ne j\le m}a_{ij}$, we have 
$$\alpha\le \sum_{i=1}^m\Vert x_i\Vert.$$ 
Put $s=2k$, ${\cal E}=\{ 
E\subset \{1,\dots,m \}:  \vert E \vert = s \}$. Write, for $E\in \cal E$, 
$$\alpha(E) = \sum_{i\in E}\sum_{j\in E\sim  \{i \}}a_{ij}. $$ It is easy 
to see that $\sum_{E\in\cal E}\alpha(E)={{m-2}\choose{s-2}}\alpha =\vert 
{\cal E}\vert{s\choose 2}{m\choose 2}^{-1}\alpha$, hence we can pick 
$E_0\in \cal E$ so that $\alpha(E_0)\le{s\choose 2}{m\choose 2}^{-1}\alpha$ 
.   Let $F = \{i\in E_0 :  \sum_{j\in E_0\sim  \{i \}} a_{ij}\ge {1\over 
k}\alpha(E_0)\}$. Then $F$ has at most $k$ elements and, clearly, each
$k$--element subset $D\subseteq E_0\sim  F$ has 
the required property.\hfill\qed

\proclaim Lemma 6. {\sl Let $U\colon l_1^k\rightarrow L_1(\mu)$ be a 
linear 
operator, $\Vert U\Vert\le 1$. Suppose there exist  $\mu$-measurable sets 
$F_1,\dots,F_k$ such that for }$i=1,\dots,k$ $$ \Vert 
1_{F_i}Ue_i\Vert\ge\delta,\quad\quad \sum_{1\le i\ne j \le k} \Vert 
1_{F_j}Ue_i \Vert\le\gamma, $$ {\sl where $0<\gamma<\delta$. Then there 
exists a linear operator $Q\colon L_1(\mu)\rightarrow l_1^k$ such that $QU 
= 
id_{l_1^k}$ and $\Vert Q\Vert\le(\delta-\gamma)^{-1}$.}   

{\noindent \bf 
Proof.} Define \ \ $W\colon L_1(\mu)\rightarrow l_1^k$\ \  by the formula\ 
\  $Wf 
= (\int fg_i\,d\mu)_{i=1}^k$,\ \  where 
$$g_i=1_{F_i}\,\sgn(Ue_i) \ {\rm for} \ i=1,\dots,k.$$ 
It is easy to check that $\Vert W\Vert \le 1$ 
and for $x\in l_1^k$ one has $\Vert WUx\Vert\ge(\delta-\gamma)\Vert 
x\Vert$. Therefore the operator $Q=(WU)^{-1}W$ has the required 
properties.  \hfill\qed 

{\noindent \bf Proof of Proposition 4.} We may assume that $\Vert 
T\Vert=1$. Apply Lemma 5 with $\eta={1\over 2}\delta$ and $x_i=T z_i$ for 
$i=1,\dots,m$. This yields a set $D\subset \{1,\dots,m\}$, with $\vert 
D\vert=k\ge{1\over 8}\delta m$, which satisfies the assertion of Lemma 5. 
Writing $\{z_i: i\in D \}=\{f_1,\dots,f_m\}$ we can define the operator $A$ 
by the formula $Ae_i=f_i$, for $i=1,\dots,k$. The existence of the operator 
$B$ follows then immediately from Lemma 6.\hfill\qed

We shall combine Propositions 1 and 4 in Theorem 8 below. Before doing 
that we would like to state a dual version of Proposition 4. Note that, if 
$\dim  X_1, \dim  Y_1 <\infty$, then 
$\gamma_{\lower3pt\hbox{$\scriptstyle 
 T$}}(U)=\gamma_{\lower3pt\hbox{$\scriptstyle T^*$}}(U^*)$. This follows 
from the principle of local  reflexivity ([LT], p.33).

\proclaim Corollary 7. {\sl Let $V\colon l_\infty^m\rightarrow X$ be an 
operator of 
norm $1$ such  that $\Vert Ve_i\Vert\ge\delta>0$ for $i=1,\dots,m$. Then 
$\gamma_{\lower3pt\hbox{$\scriptstyle V$}}(l_\infty^k)\le2\delta^{-1}$ 
for some $k\ge {1\over 8}\delta  m$.}

{\noindent \bf Proof.} Let $Z=X^*$. Pick norm one elements $z_1,\dots,z_m$ 
in $Z$ such that $z_i(Ve_i)\ge\delta$ for $i=1,\dots,m$. Using Proposition 4 
we obtain $\gamma_{\lower3pt\hbox{$\scriptstyle V^*$}}(l_1^k)\le 
2\delta^{-1}$ 
for some $k\ge {1\over 8}\delta m$. Since 
$\gamma_{\lower3pt\hbox{$\scriptstyle
V$}} (l_\infty^k)=\gamma_{\lower3pt\hbox{$\scriptstyle V^*$}}(l_1^k)$ this
completes the proof.\hfill\qed

{\ } 

\noindent {\bf The $L_1$ result}

 The main results of this section are Theorem 11 and Corollary 12 below. 
Corollary 12.A states roughly that in any good factorization of a natural 
embedding of 
$l_2^n\ $into $L_1(0,1)$ (for example, the embedding sending the unit 
vector 
basis of $l_2^n$  to the first n Rademacher functions) through an $L_1$ 
space, the operator between the two $L_1$ spaces preserves an 
$l_1^k$ space with $k$ exponential in $n$. We begin however with a 
theorem of a more general nature which is a corollary to Propositions 1 and 
4. The assumptions in both this theorem and Theorem 11 are stated in terms 
of factorization constants of an operator into an $L_1$ space through 
$L_p$ spaces via changes of densities. The relation between these constants 
and 
factorizations of natural embeddings was one of the main tools in [FJS]. We 
shall return to this relation in the proof of Corollary 12.

\proclaim Theorem 8. {\sl  Let $T\colon Z\rightarrow L_1(\mu)$ be a linear 
operator such that $T\ne 0$ and $C_{1,q}(T)<\infty$, where $1<p<q\le\infty$. 
Set} $\sigma=1-{{q^*}\over {p^*}}$, \ \  $\Delta = \Vert T\Vert^\sigma 
C_{1,q}(T)^{1-\sigma}/C_{1,p}(T)$.   \noindent {\sl Then}  
$$\gamma_{\lower3pt\hbox{$\scriptstyle T$}}(l_1^k)\le 
2(4\Delta)^{1/\sigma}\Vert 
T\Vert^{-1}\ \ \ \ {\rm for\ some}\ \  \ \ 
k\ge {1\over 8}(4\Delta)^{-1/\sigma} \Bigl ({{C_{1,p}(T)}\over {8\Vert 
T\Vert}}\Bigr)^{p^*}.$$

  {\noindent \bf Proof.} By Maurey's result [Ma] quoted in [FJS], for each 
$r\in (1,\infty]$ there is a nonnegative function $\phi_r\in L_1(\mu)$ such 
that $\int \phi_r\, d\mu = 1$ and $$ \Vert \phi_r^{-1}T \colon  Z\rightarrow 
L_r(\phi_rd\mu)\Vert = C_{1,r}(T) $$ (this uses the convention ${0\over 
0}=0$).   Set $\phi = {1\over 2}(\phi_p +\phi_q)$. Then for $r=q$ and $r=p$ 
$$\Vert \phi^{-1}T \colon  Z\rightarrow L_r(\phi\,d\mu)\Vert \le 
2^{1/r^*}C_{1,r}(T). $$ 
Consider the operator $T_1 = \phi^{-1}T \colon  
Z\rightarrow 
L_1(\phi\,d\mu)$. Since $C_{1,r}(T_1)=C_{1,r}(T)$ for $r\in(1,\infty]$, 
applying Proposition 1 to the operator $T_1$, we have ${\cal C}\le 2^{1/p^*}$. 
We can estimate for each $i$
$$\Vert 1_{F_i}\, T_1 z_i \Vert {\lower3pt\hbox{$\scriptstyle {L_q(\phi \, d 
\mu) }$}} \le 
\Vert T_1 \colon  Z\rightarrow L_q(\phi\,d\mu)\Vert \le 
2^{1/q^*}C_{1,q}(T).$$
 Hence we obtain elements $z_1,\dots,z_m$ 
in $Ball(Z)$ and sets $F_1,\dots,F_m$ so that 
$$\Vert 1_{F_i}\,T_1z_i  
\Vert_{L_1(\phi\,d\mu)}^\sigma \ge 2^{-{1\over p} -1} \bigl(2^{1\over q^*} 
C_{1,q}(T)\bigr)^{\sigma - 1} C_{1,p}(T) = 2^{-2} \Delta ^{-1} \Vert T 
\Vert^{\sigma}$$
and  $m>\bigl({{C_{1,p}(T)}\over {8\Vert T\Vert}}\bigr)^{p^*}$.  Therefore, 
we have
$$\Vert 1_{F_i}\, Tz_i \Vert _1 = \Vert 1_{F_i}\, T_1 z_i \Vert 
{\lower3pt\hbox{$\scriptstyle L_1(\phi \,d\mu)$}} \ge (4\Delta)^{-{1\over 
\sigma}} 
\Vert T \Vert.$$
 Now we simply 
apply Proposition 4 to $T$, $z_1,\dots,z_m$  and  $F_1,\dots,F_m$.\hfill\qed

We next state two corollaries to Theorem 8 concerning the dual situation.

\proclaim Corollary 9. {\sl  Let $U\colon C(K)\rightarrow X$ be a linear 
operator 
such that $0<\pi_r(U)<\infty$, where $1\le r<t<\infty$. Set $\sigma=1-{r\over 
t}$, $\Delta=\Vert U\Vert ^\sigma \pi_r(U)^{1-\sigma}/\pi_t(U)$.} {\it} Then 
$$\gamma_{\lower3pt\hbox{$\scriptstyle {U}$}}(l_\infty^k)\le 
2(4\Delta)^{1\over\sigma}\Vert U\Vert ^{-1}\ \  \ \ 
{\rm for \  some} \ \  \ \ k\ge {1\over 8}(4\Delta)^{-1/\sigma} 
\Bigl({{\pi_t(U)}\over{8\Vert U\Vert}}\Bigr)^t.$$

   {\noindent \bf Proof.} This follows easily from Theorem 8, because $C(K)^*$ 
is an $L_1$ space, $\gamma_{\lower3pt\hbox{$\scriptstyle 
U^*$}}(l_1^k)=\gamma_{\lower3pt\hbox{$\scriptstyle 
{U}$}}(l_\infty^k)$ and 
$C_{1,p}(U^*)=\pi_{p^*}(U)$ for $1<p\le\infty$ (see [R]).\hfill\qed

  \proclaim Corollary 10. {\sl If \  $U\colon l_\infty^N\rightarrow X$\  ,\  
$t>1$\  
and \ $\pi_t(U)= c\Vert U\Vert >0$\ , then \ 
$\gamma_{\lower3pt\hbox{$\scriptstyle 
{U}$}}(l_\infty^k)\le$\  $ 
2({4\over c})^{t^*}N^{t^*-1}\Vert U\Vert ^{-1}$ for some $k\ge 2^{t^*-
3}({c\over 8})^{t^*+t}N^{1-t^*}$. }   

{\noindent \bf Proof.} This follows by using 
Corollary 9 with $r=1$, because $\pi_1(U)\le N\,\Vert U\Vert$.\hfill\qed

The next theorem and Corollary 12 below are the main results of this section.

\proclaim Theorem 11.  {\sl  Let $T\colon Z\rightarrow L_1(\mu)$ be a 
bounded 
linear operator and let $1<p<\infty$. Let $Z_0\subseteq Z$. Suppose that 
$n=\dim  TZ_0<\infty$ and that} $$C_{1,p}(u)\ge c\Vert T\Vert >0, $$ {\sl for 
each finite rank operator $u\colon Z\rightarrow L_1(\mu)$ such that 
$u\vert_{Z_0}=T\vert_{Z_0}$.}\hfill \break {\sl If $c\ge 2^5$ and $\delta=(p-
1)n$, then $\gamma_{\lower3pt\hbox{$\scriptstyle T$}}(l_1^k)<5^\delta$ for 
some $k>5^{-\delta}\bigl(2^{1/\delta}\bigr)^n$.}

\proclaim Corollary 12.A.    {\sl Let $X$ be an $n$-dimensional subspace of 
$L_1$ for which $C_{1,p^*}(X)\le C\sqrt{p^*}$ for all $2\le p^* <\infty$.
 If $j$ is the inclusion map from $X$ into $L_1$ and $j=UW$ 
with $W\colon X\rightarrow Z,\  U\colon Z\rightarrow L_1,\  \Vert 
W\Vert\le 1 $ and $Z$ is an $L_1$ space, then 
$\gamma_{\lower3pt\hbox{$\scriptstyle U$}}(l_1^k) 
\le 5^{2D}$ for some $k \ge 5^{-2 D}\,2^{n/(2 D)}$,
where $D=(2^8C\Vert U\Vert)^2$.}

The proof of the corollary is very similar to the proof of Corollary 1.5 in [FJS];
there are, however, some changes.
Here is an outline of the proof:  Let  $Z_0=WX$ and define $p$ by 
$p^*={n\over D}$. 
If $n<2D$ the conclusion is obvious, so we may assume that $n\ge 2 D$.
Then $\delta =(p-1)n$ satisfies $\delta \le 2 D$, so, by Theorem 11, 
it is enough to prove that for any extension 
$V\colon Z\rightarrow L_1(\mu )$ of $U\vert _{Z_0}$,
 one has $C_{1,p}(V) \ge 2^5\Vert U\Vert$.

Now, for any extension $\tilde U \colon  Z \rightarrow X$ of $U\vert _{Z_0}$,  
we have $id_X=\tilde UW$ and hence, using a weak form of Grothendieck's
inequality, we obtain
$$n^{1/2} = \pi_{\lower3pt\hbox{$\scriptstyle 2$}}(id_X) \le 
\Vert W \Vert \pi_{\lower3pt\hbox{$\scriptstyle 2$}}(\tilde U)\le
2\,\gamma_{\lower3pt\hbox{$\scriptstyle 2$}}(\tilde U).$$
We can now apply Theorem 1.3 in [FJS] with $\beta = {1\over 2}n^{1/2}$
to get that for any extension $V \colon Z \rightarrow  L_1$ of
$U\vert _{Z_0}$ 
$$C_{1,p}(V) \ge 2^{-3} C_{1,p^*}(X)^{-1}n^{1/2}.$$ 
With the choice $p^*=n/D$ we get the desired estimate
$$C_{1,p}(V) \ge 2^{-3} C^{-1}D^{1/2} =
2^5\Vert U \Vert.\eqno{\qed}$$

Corollary 12.B strengthens Corollary 12.A in the same way that Theorem 5.1 
in [FJS] 
strengthens Corollary 1.5 in [FJS].

\proclaim Corollary 12.B.  Suppose that $X\subset L_1$, $\dim  X = 
n$, and $C_{1,p^*}(X) \le C\sqrt{p^*}$ for all $2\le p^* < \infty$.  Let  $Y$ be a 
Banach space whose dual has finite cotype $q$ constant $C_q(Y^*)$ and let 
$Q\colon Y \to 
L_1$ be an operator for which
$$Q (Ball (Y)) \supset Ball (X).$$
Let $Q=UW$ be any factorization of $Q$ through an $L_1$ space with $\Vert 
W \Vert \le 1$.  Then for some absolute constant $\eta$,
$\gamma_{\lower3pt\hbox{$\scriptstyle U$}}(l_1^k) \le 5^D$
for some
$k \ge 5^{-D}\, 2^{n/D}$,
where
$$D = \eta C^2qC_q^2(Y^*)\Vert U\Vert ^2.$$

\noindent {\bf Sketch of proof.} Follow the proof of Theorem 5.1 in [FJS] 
(with $r$ replaced by $p$) up to the place on p. 98 where it is
 proved that $C_{1,p}(U) \ge 2$. Of course, now we need and can assure that 
$C_{1,p}(\tilde U) \ge 2^5 \Vert U\Vert$ for any extension $\tilde U \colon Z 
\rightarrow
 L_1$ of the restriction of $U$ to $U^{-1}(X)$. 
Then apply Theorem 11.\hfill\qed

To prove Theorem 11 we need to introduce some notation and some 
preliminary results. Given two Banach spaces $Z$ and $W$ the space of all 
bounded linear operators between them is denoted by  $B(Z,W)$, while 
$F(Z,W)$ is the set of those $u\in B(Z,W)$ such that $\rank \  u<\infty$. By 
$\alpha$ we denote a norm on $F(Z,W)$ such that $\alpha(u)\le \Vert 
u\Vert$ if $\rank \  u = 1$.

Given a Banach space $W$ and numbers $n,\beta \ge 1$, let us denote by 
$q_{\lower3pt\hbox{$\scriptstyle W$}}(n,\beta)$ the least number $k$ such 
that, 
whenever $v \colon  W\rightarrow 
E$ is a continuous linear operator with $\rank \  v\le n$, there exists $P\in 
F(W,W)$ such that $v\,P = v$, $\Vert P\Vert\le\beta$ and $\rank  P\le k$. 
(Of 
course, we let $q_{\lower3pt\hbox{$\scriptstyle W$}}(n,\beta)=\infty$, if no 
such $k$ 
exists.)  

\proclaim Proposition 13. {\sl Let $T\in B(Z,W)$ and let $Z_0 \subseteq Z$, \ 
$\dim  T(Z_0) = n < \infty$. If $1\le \beta < \infty$ and
$q_{\lower3pt\hbox{$\scriptstyle  W$}}(n,\beta)<\infty$, 
then there exists $P\in F(W,W)$ such that $\Vert P\Vert\le\beta$, $\rank  
P\le q_{\lower3pt\hbox{$\scriptstyle W$}}(n,\beta)$ and} 
$$\alpha (PT)\ge \inf \{ 
\alpha(u):  u\in F(Z,W), 
u\vert_{Z_0} = T\vert_{Z_0}\}. $$

{\noindent \bf Proof.} Write $Y = \{u\in F(Z,W):  u\vert_{Z_0} = T\vert_{Z_0} 
\}$ and $A = \inf \{\alpha(u) :  u\in  Y \}$. By the Hahn--Banach theorem 
there is a norm one functional $\Phi$ on $(F(Z,W),\alpha)$ such that $\Phi(u) 
= A$ for each $u\in Y$.   Observe that if $S\colon W\rightarrow Z^{**}$ is the 
linear 
operator defined by $$(Sw)\,(z^*)=\Phi(z^*\otimes w), $$ then for all $u\in 
F(Z,W)$ one has $$ \Phi(u)= Tr( S u ) = Tr( u^{**} S). $$ Clearly, our 
assumption on $\alpha$ yields $\Vert S\Vert\le 1$. Moreover, for all $u\in  
Y$ one has $$ (T-u)^{**} S = 0.\eqno(2) $$ Indeed, since ${\rm Ker\,}
((T-u)^{**})  \supseteq ({\rm Ker\,}(T-u))^{\perp\perp} \supseteq
{Z_0}^{\perp\perp}$, it suffices  to verify that $SW \subseteq
Z_0^{\perp\perp}$. The latter inclusion is  obvious, because if $w\in W$,
$z_0^*\in Z_0^\perp$ then we have  $(Sw)(z_0^*) =\Phi(z_0^*\otimes w) =0$, 
since
$z_0^*$ annihilates $Z_0$.    Since $\rank u_0 = \dim  T Z_0 = n$, for some
$u_0\in  Y$, and since (2)  implies that $T^{**} S = u^{**} S$ for each $u\in
Y$,  we obtain that $\rank   T^{**} S \le n$. Hence, by the definition of
$q_{\lower3pt\hbox{$\scriptstyle  W$}}(n,\beta)$, there is a $P\in 
F(W,W)$ such that $T^{**} S P = T^{**} S$, $\Vert P \Vert\le \beta$ and 
$\rank  
P  \le q_{\lower3pt\hbox{$\scriptstyle W$}}(n,\beta)$.   Observe that, if $u$ 
is any 
element of $Y$, then $$ 
\alpha(P T)\ge\Phi(P T)= Tr(S P T) =Tr (T^{**} S P) = Tr (T^{**} S) = Tr(u^{**}S) 
= \Phi(u)=A. \eqno{\qed}$$

\proclaim Lemma 14. {\sl If $W=L_1(\mu)$ and $0<\epsilon<1$, then 
$q_{\lower3pt\hbox{$\scriptstyle W$}}(n,(1-\epsilon)^{-1})<({2\over 
\epsilon} +1)^n$.}

{\noindent \bf Proof.} Let $u\colon W\rightarrow E$ have $\rank $ $n$. 
Write 
$u=UQ_0$, where $Q_0\colon W\rightarrow W/{({\rm Ker\,}\  u)}$ is the 
quotient
map and  let $F=W/{({\rm Ker\,}\  u)}$.   Set $\beta=(1-\epsilon)^{-1}$. 
Suppose
first that for  some $k$ there exists an operator $Q\colon l_1^k\rightarrow F$
such that $\Vert  Q\Vert\le\beta$ and $Q(Ball(l_1^k))\supseteq Ball(F)$. By 
the
lifting  property of $l_1^k$ there is  $Q_1\colon l_1^k\rightarrow W$ such 
that
$\Vert Q_1  \Vert \le \Vert Q\Vert\le\beta$ and $Q=Q_0Q_1$. By the lifting
property of  $W=L_1(\mu)$ there is  $Q_2\colon W\rightarrow l_1^k$ such 
that
$\Vert Q_2  \Vert \le\Vert Q_0 \Vert  \le 1$ and $Q_0=QQ_2$. Let 
$P=Q_1Q_2$.
Then  $\Vert P\Vert\le \beta$, $\rank  P\le k$ and 
$$u=UQ_0=UQQ_2=UQ_0Q_1Q_2=uP. $$ Now the well--known volume 
argument shows that the unit sphere of $F$ contains an $\epsilon$--net 
(where $(1-\epsilon)^{-1}=\beta$) of cardinality $k<{1\over 
2}({2\over\epsilon}+1)^n$. Using this fact one easily constructs the operator 
$Q\colon l_1^k\rightarrow F$ with the two properties 
which we have used above.\hfill\qed

Theorem 11 can now be obtained by letting $\epsilon = {1\over 2}$ and 
replacing $c$ by $2^5$ in the following proposition.

   \proclaim Proposition 15. {\sl Let $T\colon Z\rightarrow L_1(\mu)$ satisfy 
the 
assumptions of Theorem 11, and let $0<\epsilon<1$. Then 
$$\gamma_{\lower3pt\hbox{$\scriptstyle T$}}(l_1^k)<4\Bigl({2\over{(1-
\epsilon)c}}\Bigr)^p \Bigl({2\over\epsilon}+1\Bigr)^{n(p-1)}\Vert T\Vert^{-1} 
$$ for some $k>4^{p-2}\bigl({{(1-\epsilon)c}\over 8}\bigr)^{pp^*} 
\bigl({2\over\epsilon}+1\bigr)^{-n(p-1)}$.}   

{\noindent \bf Proof.} Let 
$\beta=(1-\epsilon)^{-1}$, $W=L_1(\mu)$, $\alpha = C_{1,p}$. Thanks to 
Lemma 14, we can apply Proposition 13 which yields an operator $P$ on 
$L_1(\mu)$ such that $\Vert P\Vert<\beta$, $\rank \ P\le N={1\over 
2}({2\over\epsilon}+1)^n$ and $C_{1,p}(PT)\ge c\Vert T\Vert $. Clearly, 
$$C_{1,\infty}(PT)\le(\rank \ PT)\Vert PT \Vert\le N \Vert PT \Vert. $$ Let 
$q=\infty$. We can now estimate $\gamma_{\lower3pt\hbox{$\scriptstyle 
PT$}}(l_1^k)$, 
using Theorem 8. 
The resulting inequality, combined with the obvious relation 
$\gamma_{\lower3pt\hbox{$\scriptstyle T$}}(l_1^k)\le\Vert 
P\Vert\gamma_{\lower3pt\hbox{$\scriptstyle PT$}}(l_1^k)$, yields the 
desired 
estimates for 
$\gamma_{\lower3pt\hbox{$\scriptstyle T$}}(l_1^k)$ and $k$. \hfill\qed

{\ } 

\noindent{\bf The $C(K)$ result}

The main result of this section is Theorem 16 and in particular its Corollary 
20 which gives a local version of a result of Pe\l czy\'nski [Pe1] by showing 
that an 
operator from a $C(K)$ space which preserves a copy of $l_2^n$ also 
preserves a copy of $l_\infty ^k$ with $k$ an exponent of $n$.

\proclaim Theorem 16. {\sl Let $U\colon C(K)\rightarrow X$ be a bounded 
linear 
operator and let $1<t<\infty$.    Suppose that $E\subseteq C(K)$, $\dim  
E=n<\infty$ and let $$\pi_t(U\vert_E)= c\Vert U\Vert >0. $$ If $c\ge 2^5$ 
and $\alpha={n\over{t-1}}$, then $\gamma_{\lower3pt\hbox{$\scriptstyle 
{U}$}}(l_\infty^k)<5^\alpha\Vert 
U\Vert^{-1}$ for some $k>5^{-\alpha}\bigl(2^{1/\alpha}\bigr)^n$.}

For the proof we need a dual version of Lemma 14 and a proposition. The 
lemma
is a simplified version of a result which was proved in [FJS].

\proclaim Lemma 17. {\sl Let $F$ be a subspace of $C(K)$, $\dim \  F =n$. 
Let 
$0<\epsilon<1$ and $N={1\over 2}({2\over\epsilon}+1)^n$. Then there is 
$Q\colon C(K)\rightarrow C(K)$ such that $Qf=f$ for $f\in F$, $\Vert 
Q\Vert<(1-
\epsilon)^{-1}$ and $ \rank \  Q\le N$.}

{\noindent \bf Proof.} There exists $k\le N$ and an operator $J\colon 
F\rightarrow 
l_\infty^k$ such that $\Vert J\Vert<(1-\epsilon)^{-1}$ and $\Vert Jf \Vert\ge 
\Vert f \Vert$ for $f\in F$. By the extension property of $C(K)$, there is 
$J_1\colon l_\infty^k\rightarrow C(K)$ such that \ $\Vert J_1\Vert < (\Vert 
J\Vert(1-\epsilon))^{-1}$\  and \ $J_1(Jf)=f$\  for \ $f\in F$. We let\  
$Q=J_1J_2$,\ where\break $J_2\colon C(K)\rightarrow l_\infty^k$ is a linear 
extension of $J$ with $\Vert J_2\Vert =\Vert J\Vert$.\hfill\qed

Theorem 16 follows easily from the next proposition by letting $\epsilon = 
{1\over 2}$ and  replacing $c$ by $2^5$.

\proclaim Proposition 18. {\sl  Let $U\colon C(K)\rightarrow X$ satisfy the 
assumptions of Theorem 16, and let $0<\epsilon<1$. Then 
$$\gamma_{\lower3pt\hbox{$\scriptstyle 
{U}$}}(l_\infty^k)<4\Bigl({2\over{(1-
\epsilon)c}}\Bigr)^{t^*} 
\Bigl({2\over\epsilon}+1\Bigr)^{n/{(t-1)}}\Vert U\Vert^{-1} $$ for some 
$k>4^{t^*-2}\bigl({{(1-\epsilon)c}\over 8}\bigr)^{tt^*} 
\bigl({2\over\epsilon}+1\bigr)^{-n/{(t-1)}}$.}

{\noindent \bf Proof.}   By Lemma 17, there exists $P\colon C(K)\rightarrow 
C(K)$ 
such that $\Vert P \Vert<(1-\epsilon)^{-1}$, $\rank  P< N= {1\over 2} 
\bigl({2\over\epsilon}+1\bigr)^n$ and $Pe=e$ for $e\in E$. It follows that 
$$C_{1,t^*}(P^*U^*) \ge \pi_t(UP) \ge \pi_t(UP\vert_E) = c\Vert U\Vert . $$ 
Since also $C_{1,\infty}(P^*U^*)\le (\rank \ UP)\Vert (UP)^* \Vert\le N\Vert 
UP \Vert $, we can finish the proof by applying an argument similar to that 
in 
the proof of Proposition 15 and dualizing.\hfill\qed

 \proclaim Corollary 19. {\sl Let $X$ be a Banach space. Suppose that      
$\,t,n >1\,$\  and there is an operator $U\colon C(K)\rightarrow X$ and a 
subspace 
$E\subseteq C(K)$, $\dim  E=n<\infty$ such that $\pi_t(U\vert_E)= c\Vert 
U\Vert >0$, where $c\ge 2^5$.}  {\sl Write $\alpha={n\over{t-1}}$, $c_1=(\log 
2)(\log {4\over 3})/\log 5$. Then, for all 
$j\le\min\{5^{-\alpha}\bigl(2^{1/\alpha}\bigr)^n, 
{3\over 4}\exp ( {{c_1n}\over{\alpha^2} })\} $, 
$X$ contains a subspace $X_j$ such that $d(X_j,l_\infty^j)<2$.}   

{\noindent \bf Proof.} Since $c\ge 32$ using Theorem 16 we obtain that, if 
$\alpha={n\over{t-1}}$, then $\gamma_{\lower3pt\hbox{$\scriptstyle 
{U}$}}(l_\infty^k)<5^\alpha\Vert U\Vert 
^{-1}$ for some $k>5^{-\alpha}\bigl(2^{1/\alpha}\bigr)^n$. Consequently, we 
obtain the inequality $$\gamma_{\lower3pt\hbox{$\scriptstyle {id_X}$}}
(\l_\infty^k)\le \Vert 
U\Vert\gamma_{\lower3pt\hbox{$\scriptstyle {U}$}}(l_\infty^k)<5^\alpha, 
$$ from which we shall deduce a 
lower estimate for the number $$j_0=\min 
\{m: \gamma_{\lower3pt\hbox{$\scriptstyle {id_X}$}}(\l_\infty^m)\ge 2\}. 
$$ Put for 
brevity 
$g_i(X)=\gamma_{\lower3pt\hbox{$\scriptstyle {id_X}$}}(\l_\infty^i)$. We 
shall employ the estimate 
$$g_{ij}(X)\ge g_i(X){2\over{1+{g_j(X)^{-1}}}}\thinspace,\eqno(3) $$ for 
$i,j=1,2,\dots$ , which is the quantitative statement  of results of James [J] 
and 
Giesy [G] (see, e.g., [F]).

  Suppose that $j_0\le k$ and let $m$ be an integer such that 
$j_0^{m+1}>k\ge j_0^m$. Write  $A=g_{j_0}(X)\ge 2$, $B=2/(1+A^{-1})$. Using 
repeatedly (3), one obtains $$g_k(X)\ge g_{j_0}^m(X)\ge AB^{m-
1}\ge{A\over{B^2}}B^{{\log k}\over{\log j_0}}. $$ Since $A\ge 2$, we have 
${A\over{B^2}}\ge{9\over 8}$, $B\ge {4\over 3}$. Taking logarithms of both 
sides we obtain the estimate $$\alpha \log 5\ge \log(g_k(X))> (\log {4\over 
3}){{\log k}\over{\log j_0}}. $$ Using now the estimate $k>5^{-
\alpha}\bigl(2^{1/\alpha}\bigr)^n$, we get easily $$\log j_0 > (\alpha \log 
5)^{-1}(\log {4\over 3}) ({n\over\alpha}\log 2 - \alpha\log 5)= (\log {4\over 
3})({\log 2\over \log 5}{n\over{\alpha^2}}-1)=:  \log M. $$   The latter 
estimate implies that if $1<j\le \min \{k,M\}$ then $g_j(X)<2$. This implies 
that $X$ contains a subspace $X_j$ with $d(X_j,l_\infty^j)<2$ and completes 
the  proof. \hfill\qed

   \proclaim Corollary 20. {\sl Let $U\colon C(K)\rightarrow X$ be a linear 
operator 
of norm $1$. Suppose that for some subspace $E\subseteq C(K)$ such that 
$d(E,l_2^n)=a$, $\,n\ge 2\,$, one has $\Vert Ux \Vert\ge b\Vert x \Vert$ for 
$x\in E$, where $b>0$.}   {\sl Then $U$ is bounded from below by $A_1^{-
(a/b)^2}$ 
on a subspace $G\subseteq C(K)$ such that $d(G,l_\infty^j)<2$ and $j\ge 
A_2^{(b/a)^4 n}$, where $A_1,A_2$ are absolute constants $>1$.}   

{\noindent 
\bf Proof.} Since $\pi_t(l_2^n)\ge\sqrt{n\over t}$ for $t\ge 1$ [Pe2], we can 
estimate $$\pi_t(U\vert_E)\ge {b\over a}\pi_t(id_{l_2^n})\ge{b\over a} 
\sqrt{n\over t}\Vert U\Vert . $$   Assume first that $n\ge (2^5a/b)^4$. 
Letting $t=(2^{-5}b/a)^2 n$ we obtain the estimate $\pi_t(U\vert_E)\ge 
2^5\Vert U\Vert$. Let $\alpha={n\over {t-1}}$. Using Theorem 16 we obtain, 
for some $k>5^{-\alpha}2^{n/\alpha}$, a pair of operators 
$A\colon l_\infty^k\rightarrow C(K)$ and $B\colon X\rightarrow l_\infty^k$ 
such that 
$\Vert A\Vert\Vert B\Vert<5^\alpha$ and $BUA=id_{l_\infty^k}$. Let 
$F=A(l_\infty^k)$. Clearly, one has $\Vert Ux\Vert \ge 5^{-\alpha}\Vert x 
\Vert $ for $x\in F$, and $d(F,l_\infty^k)<5^\alpha$.   Now, if 
$d(F,l_\infty^k)\ge 2$, then the argument used in the proof of Corollary 19 
can be applied to $F$. Since $g_k(F)<5^\alpha$ this will produce a subspace 
$G\subseteq F$ such that $j= \dim  G > {3\over 4}\exp(c_1{n\over 
{\alpha^2}})-1$ and $d(G,l_\infty^k)< 2$. This yields the following 
conditions on numbers $A_1$, $A_2$ $$5^{n/(t-1)}\le 
A_1^{(a/b)^2},\quad\quad j\ge A_2^{(b/a)^4n}. $$   If $n<(2^5a/b)^4$, then 
we let $G$ be any $2$--dimensional subspace of $E$, so that 
$d(G,l_\infty^2)<2$. This gives the following conditions on numbers $A_1$, 
$A_2$ $$b\ge A_1^{-(a/b)^2},\quad\quad 2\ge A_2^{(2^5)^4}. $$ It is not 
difficult to check that one can find $A_1, A_2 >1$ which satisfy all the above 
conditions.\hfill\qed

{\ } 

\noindent{\bf  A space with very non-unconditional structure}

 The main result here is Theorem 21 which, roughly speaking, shows the 
existence of an $m$-dimensional space $G$ which is contained in an 
$n$-dimensional space $Z$ with an unconditional basis only if $n$ is an 
exponent 
of $m$. Moreover any such $Z$ must contain $l_\infty^k$ with $k$ an 
exponent of $m$. This solves part of problem 11.4(b) in [Pe3].

Recall that the $gl$ norm of a linear operator $T\colon X\rightarrow Y$ is 
defined 
by the formula $$gl(T)= \sup\{\gamma_{\lower3pt\hbox{$\scriptstyle 
1$}}(UT):  
U\colon Y\rightarrow l_2,\quad 
\pi_1(U)\le 1  \}, $$ and that one writes $gl(X)=gl(id_X)$. Recall [GL] also that 
the 
unconditional constant of $X$ is greater than or equal to $gl(X)$.

   \proclaim Theorem 21. {\sl There is $\delta >0$ such that for each $m\ge 
2$ there is a Banach space $G_m$, $\dim  G_m =m$, with the following 
property. If $Z$ is a Banach space which contains an iso\-metric copy of 
$G_m$, 
then there is a subspace $Z_1$ of $Z$ such that 
$d(Z_1,l_\infty^k)<2$, where $k\ge \exp(\delta m\,gl(Z)^{-4})$.}

In fact, a somewhat stronger version of this theorem follows by applying 
Lemma 23 to the space obtained in Lemma 22. A stronger version of Lemma 
22 appears as Th. 7.1 in [P].   
\proclaim Lemma 22. {\sl There is a 
constant $B<\infty$ such that for $n=1,2,\dots$ there is a Banach space 
$F_n$, $\dim \  F_n =2n$, and a linear operator $v\colon l_2^n\rightarrow 
F_n$ such 
that $\Vert ve \Vert\ge\Vert e \Vert$  for $e\in l_2^n$ and $\pi_1(v^*)\le 
B$.}  

\noindent {\bf Proof.}  Consider a linear isometry $U\colon 
L_2^{3n}\rightarrow 
L_2^{3n}$. Write $$E_1=U([e_1,\dots,e_{2n}]), \quad\quad 
E_2=U([e_{n+1},\dots,e_{3n}]). $$ It is well--known (see e.g. [K], [S], or [P] 
Cor.7.4) that for ``most choices" of $U$ one has $$\Vert f \Vert_{L_2^{3n}}\le 
b \Vert f \Vert_{L_1^{3n}}, $$ for $f\in E_1\cup E_2$, where $b$ can be 
taken independent of $n$.   Let us fix a pair $E_1,E_2$ with the latter 
property. Put $F_n=L_1^{3n}/E_2^\perp$ and let 
$$u\colon E_1/E_2^\perp\rightarrow L_1^{3n}/E_2^\perp = F_n $$ be the 
natural 
map. If $E_1/E_2^\perp$ is given the norm induced from 
$L_2^{3n}/E_2^\perp$, then $E_1/E_2^\perp$ is isometric to $l_2^n$ and our 
choice of $E_1$ yields the estimate $\Vert e \Vert\le b \Vert ue \Vert$ for 
$e\in E_1/E_2^\perp$.   Now $u^*$ can be regarded as the composition of the 
embedding map $i_{\infty,2}^{E_2}\colon (E_2)_\infty\rightarrow(E_2)_2$ 
with the 
orthogonal projection $P$ from $L_2^{3n}$ onto $E_1\cap E_2$. Hence our 
choice of $E_2$ yields 
$$\pi_1(u^*)\le \Vert P\Vert 
\pi_1(i_{\infty,2}^{E_2})\le \pi_1(i_{\infty,1}^{E_2}) \Vert i_{1,2}^{E_2} \Vert 
\le b. $$ This shows that the operator $v=bu$ has the required properties, if 
$B=b^2$.\hfill\qed

   \proclaim Lemma 23. {\sl  Let $F$ be a Banach space and let 
$v\colon l_2^n\rightarrow F$, $1<t<\infty$. Suppose that $$\pi_t(v)\ge 
c\Vert v 
\Vert >0,\quad\quad \pi_1(v^*)\le B\Vert v \Vert. $$ Let $j\colon 
F\rightarrow Z$ 
be a linear operator such that $\Vert jf \Vert\ge\Vert f \Vert$  for $f\in 
v(l_2^n)$.}   {\sl If $c\ge 2^5B\,gl(j)$, then $Z \supset Z_1$  such that
$d(Z_1,l_\infty^k)<2$ and $k=\dim  Z_1\ge 
\min\{5^{-\alpha}\bigl(2^{1/\alpha}\bigr)^n, {3\over 4}\exp ( 
{{c_1n}\over{\alpha^2} })\} $, where $\alpha$ and $c_1$ are as in Corollary  
19.}

 \noindent {\bf Proof.} Observe that, since $gl(j^*)=gl(j)$, one has 
$$\gamma_{\lower3pt\hbox{$\scriptstyle 
{\infty}$}}(jv)=\gamma_{\lower3pt\hbox
{$\scriptstyle 1$}}((jv)^*)\le\pi_1(v^*)gl(j^*)\le B\,gl(j)\Vert 
v\Vert. $$ Consider a  \ $C(K)$--factorization of \  $jv\colon 
l_2^n\rightarrow 
Z^{**}$\ , say \ $jv=Ui$\ , where  $$\Vert i\colon l_2^n\rightarrow 
C(K)\Vert=1, \qquad \Vert U\colon C(K)\rightarrow Z^{**}\Vert\le 
B\,gl(j)\Vert 
v\Vert.$$'
Put  $E=i(l_2^n)$, then $$\pi_t(U\vert_E)\ge \pi_t(Ui)\ge \pi_t(v)\ge c\Vert 
v 
\Vert\ge 2^5\Vert U\Vert . $$ Since 
$\gamma_{\lower3pt\hbox{$\scriptstyle 
{id_{Z^{**}}}$}}(l_\infty^k)=\gamma_
{\lower3pt\hbox{$\scriptstyle  {id_{Z}}$}}(l_\infty^k)$, the conclusion 
follows from Corollary 19.\hfill\qed

  \noindent {\bf Proof of Theorem 21.}  We may assume that 
$m>2(2^5B\,gl(Z))^4$, 
where $B$ is the constant from Lemma 22. (If not, we just let 
$\delta={1\over 4} (2^5B)^{-4}$ and $G_m$ can be arbitrary space of 
dimension $m$.)   Let $G_{2n}=F_n$ and $G_{2n+1}=F_n\oplus l_1^1$ for 
$n\ge(2^5B)^4$. Fix an $m$ and let $Z$ be a Banach space and 
$j\colon G_m\rightarrow Z$ an isometric embedding, so that $gl(j)\le \,gl(Z)$.
 Let $v\colon l_2^n\rightarrow G_m$ be the operator 
obtained from that in Lemma 22 (here $m=2n$ or $m=2n+1$). Observe that, 
for $t\ge 1$, $$\pi_t(v)\ge\pi_t(id_{l_2^n})\ge\sqrt{n\over t}. $$ Letting 
$t=(2^5B\,gl(j))^{-2}n$ and applying Lemma 23, one can easily find the 
absolute 
constant $\delta$ needed in Theorem 21.\hfill\qed

\noindent {\bf Remark.} It follows from ([P], Th. 7.1) that the spaces $G_n$ 
in Theorem 21 
can be chosen to have uniform cotype 2 constants.

 \beginsection References

\item{[F]} T. Figiel, {\it Factorization of compact operators and  
applications to the approximation problem,} {\bf Studia Math. 45}  
(1973), 191-210.  

\item{[FJ]} T. Figiel and W. B. Johnson, {\it Large subspaces of  
$l_{\infty}^n$ and estimates of the Gordon-Lewis constant,} {\bf Israel J. 
Math. 37} (1980), 92--112.

\item{[FJS]}  T. Figiel, W.~B. Johnson, and G. Schechtman, {\it  
Factorizations of natural embeddings of  $l_p^n$ into  $L_r$, I,} {\bf  
Studia Math. 89} (1988), 79--103. 

 \item{[G]}  D. P. Giesy,  {\it On a convexity condition in normed linear 
spaces,} 
{\bf Trans. Amer. Math. Soc. 125}  (1966), 114-146.

\item{[GL]} Y.~Gordon and D.~R.~Lewis, {\it Absolutely summing operators 
and local 
unconditional structures,} {\bf Acta Math. 133} (1974), 27-48.

\item{[J]} R. C. James, {\it Uniformly non-square Banach spaces,} {\bf 
Ann. of Math. 80} (1964), 542--550.    

\item{[JS]} W. B. Johnson and G. Schechtman, {\it On subspaces of  $L_1$ 
with maximal distances to Euclidean space,} {\bf  Proc. Research 
Workshop on Banach Space Theory, Univ. of Iowa}, Bor-Luh Lin, ed. 
(1981), 83--96.   

\item{[K]} B. S. Kashin, {\it Sections of some finite-dimensional sets and 
classes of smooth functions,} {\bf  Izv. Akad. Nauk SSSR 41} (1977), 
334--351 (Russian).  

\item{[LT]} J. Lindenstrauss and L. Tzafriri, {\it Classical Banach spaces I, 
Sequence spaces,}  Springer-Verlag, (1977). 

\item{[M]} B. Maurey, {\it Th\'eor\`emes de factorisation pour les 
op\'erateurs lin\'eaires \`a valeurs dans les espaces  $L_p$,} {\bf  
Asterique No. 11, Soc. Math. France} (1974).   

\item{[P]} G. Pisier, {\it Factorization of linear operators and geometry of 
Banach spaces,} {\bf  CBMS Regional Conference Series in Mathematics 60} 
(1986), AMS, Providence. 

\item{[Pe1]} A. Pe\l czy\'nski, {\it Projections in certain Banach spaces,} {\bf  
Studia 
Math. 19} (1960), 209--228.

\item{[Pe2]} A. Pe\l czy\'nski, {\it A characterization of Hilbert-Schmidt 
operators,} {\bf  
Studia Math. 58} (1967), 355--360.

\item{[Pe3]} A. Pe\l czy\'nski, {\it Finite dimensional Banach spaces and 
operator ideals,} 
in {\bf Notes in Banach Spaces}, H.~E. Lacey Ed., Univ. of Texas (1980), 81--
181.

\item{[R]}  H. P. Rosenthal, {\it On subspaces of  $L^p$,} {\bf Ann. of 
Math. 97} (1973), 344--373. 

\item{[Sz]} S. J. Szarek, {\it On Kashin's almost Euclidean orthogonal 
decomposition of  $l_1^n$,} {\bf Bull. Acad. Polon. Sci. 26} (1978), 691--694.

{\baselineskip=14pt 

{\  }

{\ }

 {\obeylines \parindent=0pt \lineskip=2pt
Institute of Mathematics, 
Polish Academy of Sciences, 
Abrahama 18, 81-825 Sopot, Poland

{\ }

Department of Mathematics, 
Texas A\&M University,
College Station, Texas 77843, USA  

{\ }

Department of Theoretical Mathematics, 
The Weizmann Institute of Science, 
Rehovot, Israel 

{\ } }}

\vfil
\end